\documentclass[letterpaper, 10 pt, conference]{ieeeconf}  

\IEEEoverridecommandlockouts                              

\usepackage{amsmath, amsfonts}
\usepackage{color}
\usepackage{amssymb}
\usepackage{mathtools}
\usepackage{csquotes}
\usepackage{quotes}
\usepackage{float}

\usepackage{balance}
				
\newtheorem{theorem}{Theorem}

\usepackage{rotating}
\usepackage{adjustbox}

\usepackage{calc}

\usepackage{longtable}
\usepackage{multirow}

\usepackage{cite}

\usepackage{graphicx}

\usepackage{placeins}

\usepackage{tabularx}

\usepackage{lmodern,textcomp}

\usepackage{enumerate}

\usepackage{xspace}

\newcommand*\tageq{\refstepcounter{equation}\tag{\theequation}}

\usepackage{algorithm}
\usepackage{algpseudocode}

\newcommand{\blambda}{\boldsymbol{\lambda}\xspace}

\renewcommand{\t}{^\textsf{T}\xspace}

\newcommand{\away}[1]{}

\newcommand{\R}{\mathbb{R}\xspace}

\newcommand{\N}{\mathbb{N}\xspace}

\newcommand{\cL}{\mathcal{L}\xspace}
\newcommand{\cX}{\mathcal{X}\xspace}

\newcommand{\cT}{\mathcal{T}\xspace}
\newcommand{\cO}{\mathcal{O}\xspace}

\newcommand{\bD}{\textbf{D}\xspace}

\newcommand{\bx}{\textbf{x}\xspace}

\newcommand{\by}{\textbf{y}\xspace}

\newcommand{\bg}{\textbf{g}\xspace}

\newcommand{\bz}{\textbf{z}\xspace}

\newcommand{\bei}[1]{{\textbf{e}}\xspace}

\newcommand{\bO}{\textbf{0}\xspace}

\newcommand{\opdiag}{\operatorname{diag}\xspace}

\newcommand{\commentout}[1]{}

\usepackage{datetime,url}

\usepackage{cite}

\newcommand{\abbFig}{Figure}
\newcommand{\abbfig}{Figure}

\title{\LARGE \bf
An Integral Penalty-Barrier Direct Transcription Method for Optimal Control
}

\author{Martin P. Neuenhofen$^{1}$ and Eric C. Kerrigan$^{2}$
\thanks{$^{1}$Martin P. Neuenhofen is with the Deparment of Electrical \& Electronic Engineering, Imperial College London, SW7 2AZ London, UK, e-mail: m.neuenhofen19@imperial.ac.uk
        {\tt\small www.MartinNeuenhofen.de}}%
\thanks{$^{2}$Eric C. Kerrigan is with the Deparment of Electrical \& Electronic Engineering and Department of Aeronautics, Imperial College London, SW7 2AZ London, UK, e-mail: e.kerrigan@imperial.ac.uk
        {\tt\small www.imperial.ac.uk/people/e.kerrigan}}%
}

\begin{document}

\sloppy
\maketitle
\thispagestyle{empty}
\pagestyle{empty}

\begin{abstract}
Some direct transcription methods can fail to converge, e.g.\ when there are singular arcs.
We recently introduced a convergent direct transcription method for optimal control problems, called the penalty-barrier finite element method (PBF). PBF converges under very weak assumptions on the problem instance. PBF avoids the ringing between collocation points, for example, by avoiding  collocation entirely. Instead, equality path constraint residuals are forced to zero everywhere by an integral quadratic penalty term.
We highlight conceptual differences between collocation- and penalty-type direct transcription methods. Theoretical convergence results for both types of methods are reviewed and compared. Formulas for implementing PBF are presented, with details on the formulation as a nonlinear program (NLP), sparsity and solution. Numerical experiments compare PBF against several collocation methods with regard to robustness, accuracy, sparsity and computational cost. We show that the computational cost, sparsity and construction of the NLP functions are roughly the same as for orthogonal collocation methods of the same degree and mesh. As an advantage, PBF converges in cases where collocation methods fail. PBF also allows one to trade off computational cost, optimality and  violation of differential and other equality equations against each other.
\end{abstract}


\section{Motivation}
Direct transcription is a popular method for the numerical solution of optimal control, estimation and system identification problems. Many direct transcription methods are based on 
collocation methods~\cite{Betts2nd}. It is known that the latter can struggle with singular arcs and high-index differential algebraic equalities (DAEs) \cite[Chap.~2~\&~4.14]{Betts2nd}; the former arising in the example:
\begin{align*}
\min_{y,u}& & J&=\int_0^{\frac{\pi}{2}} \left(\,y(t)^2 + \cos(t)  u(t)\,\right)\,\mathrm{d}t,\tag{X1}\label{eqn:ExampleOCP}\\
\text{s.t.}& & y(0)&=0,\ \, \dot{y}(t)=\frac{1}{2}  y(t)^2+u(t),\ \, -1.5\leq u(t)\leq 1\,,
\end{align*}
which has the optimal solution \mbox{$u^\star(t)=\frac{1}{2}-\frac{1.5}{\left(\cos(t)-2\right)^{2}}$}. 

Quadratic integral penalty methods \cite{Balakrishnan68,Hager90,PBF} are an alternative to collocation methods, where the squared path equality constraint residual is integrated and added as a penalty into the objective. In \cite{PBF} we present such a method, called the penalty-barrier finite element method (PBF), with a rigorous proof of convergence under very mild and easily verifiable assumptions, including convergence for singular arcs and high-index DAEs. The assumptions are discussed in Section~\ref{sec:PBF_conv} below.

Convergence proofs with enforceable assumptions are of high practical relevance. If convergence can only be guaranteed under unverifiable or unenforceable assumptions, then a practitioner cannot determine whether features in numerical solutions are due to optimality or  numerical failure.

\begin{figure}
	\centering
	\includegraphics[width=0.8\columnwidth]{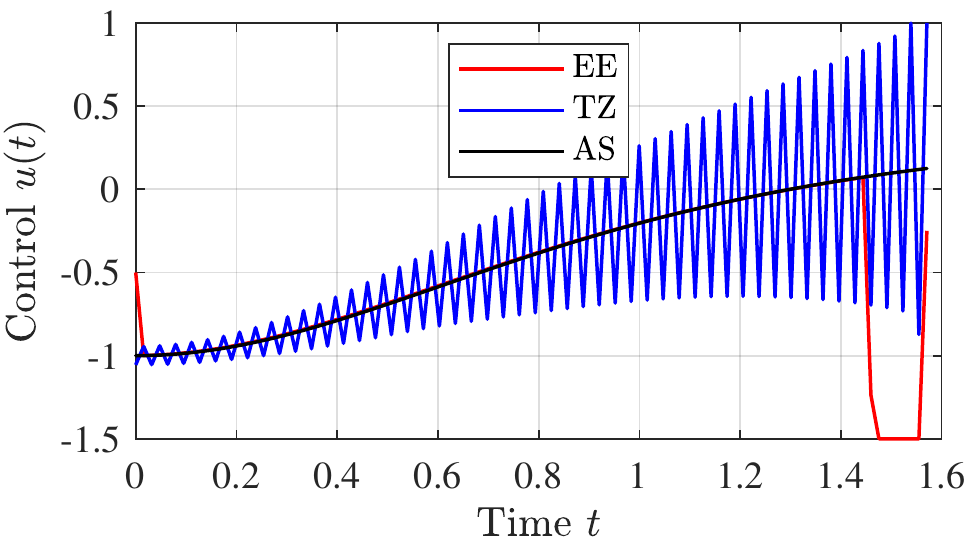}
	\caption{Solutions to \eqref{eqn:ExampleOCP}: Explicit Euler (EE), Trapezoidal (TZ), analytic solution (AS).}\label{fig:ConvExample1a}
\end{figure}

\begin{figure}
	\centering
	\includegraphics[width=0.8\columnwidth]{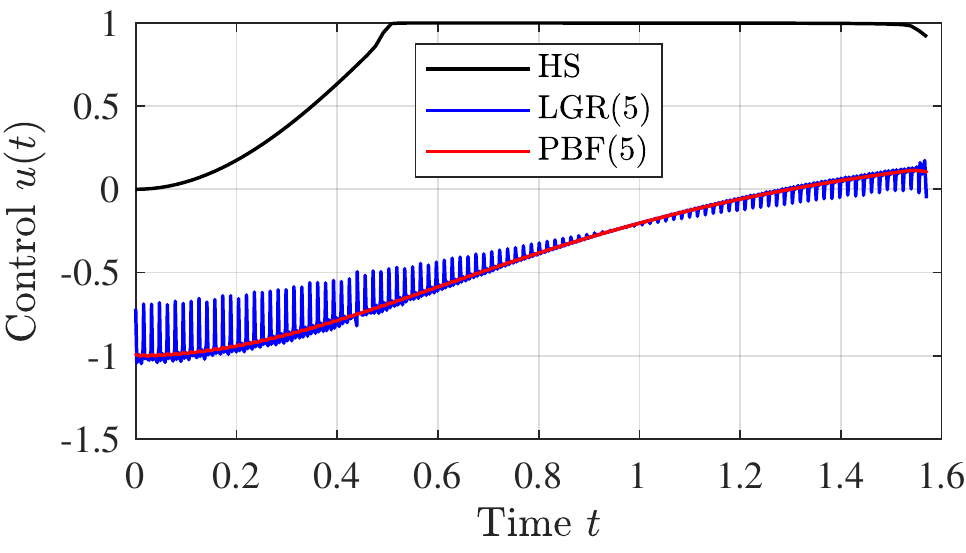}
	\caption{Solutions to \eqref{eqn:ExampleOCP}: Hermite-Simpson (HS), Gauss-Legendre-Radau of degree $5$ (LGR5), penalty-barrier finite element method of degree $5$ (PBF5).}\label{fig:ConvExample1b}
\end{figure}

Figures~\ref{fig:ConvExample1a}--\ref{fig:ConvExample1b} illustrate this dilemma by showing numerical solutions to \eqref{eqn:ExampleOCP}. They show solutions for the explicit Euler method (EE), trapezoidal method~(TZ), Hermite-Simpson collocation~(HS), Legendre-Gauss-Radau collocation of degree 5~(LGR5), and PBF of degree 5~(PBF5), each using 100 mesh intervals of equal length.
PBF5 converges, while EE diverges for the last section of the time interval. The box constraints have been imposed only to avoid unboundedness of the NLP resulting from the discretization. The solutions from HS and LGR2 (the latter of which was not plotted) also diverge if box constraints are removed. For LGR5 we see the best result among the non-penalty-type methods, however the solution suffers from ringing. All NLPs have been solved until the unscaled KKT residual $\infty$-norm of $<10^{-10}$ stagnated.

\subsection{Contributions}
In this paper we summarise the presentation in \cite{PBF} and extend it with practical aspects and new examples.

First, we generalize the problem statement. We extend the objective with a Mayer term and incorporate assumptions on the problem such that the convergence analysis from \cite{PBF} applies. We extend the box-constraints to both the states and controls. We explain how these constraints can be posed into the NLP via a projection matrix, and where this matrix arises in the optimality conditions.

We then draw a connection between the first-order optimality conditions between the PBF discretization and the first-order optimality conditions of a regularized primal-dual problem.
In this context, we visualize and compare sparsity patterns of the derivatives in PBF with the LGR direct collocation method.

We cover further practical aspects: We illustrate the stability condition for quadrature, which is a building block when constructing a convergent implementation of PBF. In the numerical experiments, we demonstrate how PBF minimizes a bias between optimality and feasibility with respect to the choice of the penalty parameter.

\subsection{Outline}
Section~\ref{sec:NOC} reviews the problem statement and notation.
Section~\ref{sec:CvP} compares collocation-type methods with integral penalty-type methods. We highlight key differences in the numerical approaches and review available convergence results.
Section~\ref{sec:NLP} shows how PBF can be implemented, with emphasis on the formulation of the NLP and the sparsity of derivative matrices.
Section~\ref{sec:NumExp} presents numerical experiments, where we compare several collocation methods with PBF and focus on numerical cost, sparsity and (order of) convergence.

\section{Numerical Optimal Control}\label{sec:NOC}

\subsection{Problem Statement}
A  common reference problem format for optimal control is the Bolza form~\cite{KellyMatthew}, which can cover most optimal control, estimation and system identification problem instances, including multi-phase problems~\cite{GPOPS2}. Without loss of generality, the initial and end time can be fixed~\cite{Betts2nd}. 
Using slack variables, constraints can be separated into nonlinear equality constraints and box inequality constraints, resulting in the form
\newcommand{\xL}{{x_\textsf{L}}\xspace}
\newcommand{\xR}{{x_\textsf{U}}\xspace}
\newcommand{\refOCP}{(OCP)\xspace}
\begin{subequations}
	\begin{align}
	\operatornamewithlimits{min}_{x:=(y,u)\in\cX}&&W\big(y(t_0),{}&{}y(t_F)\big)+\int_\Omega w\big(y(t),u(t),t\big)\,\mathrm{d}t\tag{OCPa}\label{eqn:OCPa}\\
	\text{s.t.}& 	&0= {}&{}c\big(\dot{y}(t),y(t),u(t),t\big)\quad\text{f.a.e. }t\in\Omega\,,\tag{OCPb}\label{eqn:OCPb}\\
	& 				&0= {}&{}g\big(y(t_0),y(t_F)\big)\,,\tag{OCPc}\label{eqn:OCPc}\\
	& 				&\xL(t)\leq{}&{} x(t)\leq \xR(t)\quad\text{f.a.e. }t\in\Omega\,,\tag{OCPd}\label{eqn:OCPd}
	\end{align}
\end{subequations}
where the constraints have to be satisfied for almost every $t\in \Omega:=[t_0,t_F]\subset\mathbb{R}$.  Here \eqref{eqn:OCPd} is optional, or can be specified as only upper/lower bounds or only over subintervals of $\Omega:=[t_0,t_F]$ for each  component of~$x$.

The trajectory $x:=(y,u)$ consists of continuous functions $y : \Omega \rightarrow \R^{n_y}$ and (potentially) discontinuous functions $u : \Omega \rightarrow \R^{n_u}$. Many works \cite{RaoHager,Conway2012,KellyMatthew} separate \eqref{eqn:OCPb} into ODEs and algebraic constraints of the form
\begin{align*}
\dot{y}(t)={}&{} f_1\big(y(t),u(t),t\big)\,,\tag{OCPb1}\label{eqn:OCPb1}\\
 		 0={}&{} f_2\big(y(t),u(t),t\big)\,,\tag{OCPb2}\label{eqn:OCPb2}
\end{align*}
with functions $f_1,f_2$, describing the dynamic constraints.
This is necessary if the ODE is solved with a different scheme or when the convergence analysis poses different criteria on $f_1$ and $f_2$.

In either case, it can be considered an advantage of PBF to be capable of solving \eqref{eqn:OCPb} directly, i.e.\ without algebraic transformation into \eqref{eqn:OCPb1}--\eqref{eqn:OCPb2}, because this allows avoiding divisions by zero and phrasing equations such that their derivatives are smoother, which improves the efficiency of an NLP solver.

\subsection{Notion of Convergence}\label{sec:NotionOfConv}

Numerical methods must minimize
\begin{align*}
J(x) := &W\big(y(t_0),y(t_F)\big)+\int_\Omega w\big(y(t),u(t),t\big)\,\mathrm{d}t
\end{align*}
subject to forcing the total constraint violation
\begin{align*}
r(x) := \int_\Omega \left\|c\big(\dot{y}(t),y(t),u(t),t\big)\right\|_2^2\,\mathrm{d}t + \left\|g\big(y(t_0),y(t_F)\big)\right\|_2^2
\end{align*}
to zero while obeying \eqref{eqn:OCPd}. To keep their computation time finite, all methods search for an approximation to the optimal $x$ in a finite-dimensional space $\cX_{h,p}$, which is a subspace of the solution space $\cX$ of continuous functions~$y$ and discontinuous functions $u$; here $h$ is the mesh step size and $p$ is the degree of a method, in the case where the method is parametric. Formal details on the solution spaces are given in \cite{PBF}.

When solving optimization problems, implying \refOCP, there are two measures to take into account: feasibility and optimality. One wants to find a weak numerical solution $x_h \in \cX_{h,p}$ that solves \refOCP in a \textit{tolerance-accurate} sense. This means, the \textit{optimality gap}
\begin{align}
g_\text{opt} :=& \max\lbrace 0,J(x_h)-J(x^\star)\rbrace\,,
\end{align}
where $x^\star$ is a global minimizer of \refOCP, and \textit{feasibility residual}
\begin{align}
r_\text{feas} :=& r(x_h) \label{eqn:residual}
\end{align}
converge to zero as the mesh step size $h$ decreases. Order-of-convergence results that specify the rate at which $g_\text{opt}$ and~$r_\text{feas}$ approach zero in terms of $h$ can be derived; see~\cite{PBF}.

\section{Collocation vs Integral Penalty Methods}\label{sec:CvP}
Many state-of-the-art methods for the solution of~\refOCP are collocation-type methods.
We describe their geometric approach below and discuss their principle issues.
We then introduce integral penalty methods and explain how they can overcome the convergence problems that collocation methods can exhibit.
The section ends with a direct comparison of theoretical convergence results.

\subsection{Collocation-type Methods}
These methods determine solutions to $y$ by solving~\eqref{eqn:OCPb1} at a finite number of time instances. The trajectory $u$ must then be found such that $(y,u)$ satisfies~\eqref{eqn:OCPb2}. Examples of such methods (see~\cite{rao_asurvey} for a survey) are: single- and multiple-shooting methods \cite{Quirynen2017} based on 
certain classes of Runge-Kutta or other collocation methods for solving the differential equations; multi-step methods, in particular backward-differentiation formula methods; and direct collocation methods, an important class being orthogonal direct collocation methods.

All these methods replace the uncountable set of path constraints \eqref{eqn:OCPb1}--\eqref{eqn:OCPb2} by a finite set of constraints, namely
\newcommand{\bbThp}{\mathbb{T}_{h,p}\xspace}
\newcommand{\tbbThp}{\tilde{\mathbb{T}}_{h,p}\xspace}
\begin{align*}
\dot{y}(t)&=f_1\left(y(t),u(t),t\right) & \forall t &\in \tbbThp,\\
0&=f_2\left(y(t),u(t),t\right) &  \forall t &\in \bbThp,
\end{align*} 
where $\tbbThp,\bbThp$ are finite subsets of the interval $[t_0,t_F]$. In pure collocation methods $\tbbThp=\bbThp$, and these points are called collocation points. In some Runge-Kutta and linear multi-step methods, $\tbbThp$ is determined from the step size and Butcher tableau coefficients.

Depending on the method, various finite-difference expressions are used for the derivative $\dot{y}$ on $\tbbThp$. Collocation methods (including certain implicit and explicit Runge-Kutta methods), and linear multi-step methods use derivatives of piecewise polynomial interpolants. General Runge-Kutta methods (which are nonlinear) may use nonlinear interpolants.

Since the differential equation is satisfied only at a finite number of points and~$\cX_{h,p}$ is a strict subset of~$\cX$, the feasibility residual $r_\text{feas}$ is non-zero, in general. If the residual is too large, then the dimension of $\cX_{h,p}$ is usually increased in correspondence with choosing a set $\mathbb{T}_{h,p}$ with more elements (as in $h$-methods) and/or increasing the degree of the polynomials in $\cX_{h,p}$ (as in $p$- or $hp$-methods).
It is well-known that if care is not taken with the choice of $\mathbb{T}_{h,p}$ and the degree of the polynomials in $\cX_{h,p}$, then the feasibility residual will not converge, e.g.\ due to Runge's phenomenon. It is also possible that a solution might not exist, e.g.\ when there are more constraints than degrees of freedom.

Essentially, problems arise in collocation-type methods because they cannot explicitly take into account what happens in-between collocation points. Problems show up in two ways: On the one hand, convergence analyses rely on very technical assumptions in order for the numerical solution to converge. On the other hand, numerical solutions ring or diverge when these assumptions are violated.

Finally, collocation discretizations yield infeasible NLPs when the OCP is consistently overdetermined. Consider, for example
\begin{align*}
	y(0)=1\,,\quad\dot{y}(t)=u(t)\,,\quad\bO=\begin{bmatrix}
	\exp(t)-u(t)\\
	y(t)-u(t)
	\end{bmatrix}.
\end{align*}
This set of constraints has the consistent analytic solution $y^\star(t)=\exp(t)$. However, using any collocation-type method yields $u_h(t)=\exp(t)$ $\forall t \in \bbThp$ from the second path constraint. Further, $\dot{y}_h(t)=u_h(t)$ $\forall t \in \tbbThp$ from the ODE constraint. One would need to construct $\tbbThp$ so that $y_h(t)=\exp(t)$ $\forall t \in \bbThp$ holds, because this would be required in order for the second path constraint to be satisfied. In the case of pure collocation methods, i.e.\  $\tbbThp=\bbThp$, the discretization is infeasible, because there is no exact linear quadrature rule for the exponential.

\subsection{Integral Penalty-type Methods}
Integral penalty-type methods differ fundamentally from collocation-type methods by explicitly considering the violation of constraints everywhere in $[t_0,t_F]$ in the formulation of the optimization problem. This is possible due to the observation that one can replace the uncountable set of constraints~\eqref{eqn:OCPb} with the finite set of constraints
\begin{equation}
\int_{T} \|c\left(\dot{y}(t),y(t),u(t),t\right)\|_2^2 \mathrm{d}t \leq \varepsilon_T,\ \forall T \in \cT_h\,.
\tag{OCPb'}\label{eqn:daeint}
\end{equation}
Here $\cT_h$ is a finite set of disjoint intervals such that $\bigcup_{T\in\cT_h} \operatorname{cl}(T) = [t_0,t_F]$. Clearly,~\eqref{eqn:daeint} is satisfied with $\varepsilon_T=0$ only if~\eqref{eqn:OCPb} is satisfied; however, since~$\cX_{h,p}$ is finite dimensional, usually the best that one can do is to satisfy the above constraint with  $\varepsilon_T>0$. This also acknowledges the fact that one cannot assert feasibility of a discretization when forcing the constraint residual exactly to zero at any particular point.

The approach of integral penalty-type methods is therefore to augment the objective functional with a weighted version of the integrals in~\eqref{eqn:daeint} and a quadratic penalty on the equality constraints~\eqref{eqn:OCPc}, namely~$r$ as defined in~\eqref{eqn:residual}.

The inequality constraints~\eqref{eqn:OCPd} also need careful consideration. Here collocation methods again only enforce these at a finite number of points $\mathbb{T}_{h,p}$. This can result in overshooting of the numerical solution in-between collocation points.

Our penalty-barrier method in \cite{PBF} treats \eqref{eqn:OCPd} by expressing them as positivity constraints, introducing suitable slack variables $s : \Omega \rightarrow \R^{n_s}_+$, and then augmenting the objective functional with the \emph{integral} logarithmic barrier function
\begin{align}
\Gamma(x):=&-\sum_{j=1}^{n_s} \int_\Omega \log\big( s(t) \big)\,\mathrm{d}t\label{eqn:def:Gamma}
\end{align}
to ensure that the inequality constraints are satisfied over the whole interval $\Omega$. Details on the exact reformulation with the slacks and ensuring feasibility in the context of interior points are given in \cite{PBF}.

The resulting method consists of solving the penalty-barrier problem
\begin{align}
	\operatornamewithlimits{min}_{x \in \cX_{h,p}}\  J_{\omega,\tau}(x) = J(x) + \frac{1}{2\omega}r(x) + \tau \Gamma(x)\tag{PBP}\label{eqn:PBF}
\end{align}
for suitable penalty- and barrier parameters $\omega,\tau>0$ with respect to the mesh size $h>0$ and degree $p \in \N$ of the finite elements.

\subsection{Convergence Results}
We briefly compare collocation-type methods to integral penalty-type methods in terms of convergence guarantees.
Most important is whether a method converges. Desirably, this should be achieved under mild assumptions.
Thereafter, the order of convergence is of importance, so that highly accurate solutions can be obtained with an affordable amount of computing resources.

\subsubsection{Collocation Methods}
A convergence proof for a discretization based on the explicit Euler method is given in \cite{Maurer}. The result  makes the following assumptions: (i)~functions defining the problem must be locally differentiable with Lipschitz continuous derivatives; (ii)~there must be a local solution where the trajectories of the state and free variables are continuously differentiable and continuous, respectively; (iii)~a  homogeneity condition on active constraints; (iv)~surjectivity of linearized equality constraints; and (v)~a  coercivity assumption.

These conditions are sophisticated, difficult to understand, and very hard to verify and ensure by construction. 
The conditions ensure that the first order optimality conditions of the non-discretized and of the Euler-discretized optimization problem have a convergent solution. This is why convergence proofs for other collocation schemes make similar assumptions.

Since Euler's method converges only to first order, convergence results based on Euler's method are of limited practical use compared to results for higher-order methods. Convergence of a higher-order method is shown in~\cite{kameswaran_biegler_2008} by making additional assumptions: (vi)~functions $W,w,c,g$ must be sufficiently smooth; (vii)~the state and co-state trajectories must be sufficiently smooth; (viii)~the NLP arising from the discretization must satisfy the Linear Independence Constraint Qualification and Second-Order Sufficient Condition.

\subsubsection{The Penalty-Barrier Finite Element Method}\label{sec:PBF_conv}

In contrast to the above results for collocation methods, for PBF as in \cite{PBF} convergence follows under the following assumptions:
\begin{enumerate}[\text{(A.}1\text{)}]
	\item 
	\refOCP has a global minimizer $x^\star$.
	\item 
	$|w(\dot{y}(t),y(t),u(t),t)|$, $ \|c(\dot{y}(t),y(t),u(t),t)\|_1$,  
		$\|W(y(t_0),y(t_F))\|_1$, $\|g(y(t_0),y(t_F))\|_1$ are bounded for all arguments $x \in \cX$.
	\item 
	$W,\ w,\ c,\ g$ are globally Lipschitz continuous in all arguments except $t$.
	\item 
	The problem \eqref{eqn:PBF} has a bounded solution $x^\star_{\omega,\tau}$.
	\item The solution $x^\star_{\omega,\tau}$ can be approximated by a function $x_h \in \cX_{h,p}$ with an error that converges to zero as $h\searrow+0$.
\end{enumerate}
We review the discussion of these assumptions from \cite{PBF}:
\begin{enumerate}[\text{(A.}1\text{)}]
	\item is reasonable.
	
	\item can be enforced by construction. To this end, $W,w,c,g$ can be bounded below and/or above, if necessary, with minimum and maximum terms. For example, if $|J(x^\star)|\ll 10^{8}$ is expected, then we can replace $w$ by a modified function
		$$ \tilde{w}(\chi,\upsilon,t):=\max\Big\{-10^{8}\,,\,\min\big\{ w(\chi,\upsilon,t)\,,\,10^{8}\big\}\Big\}\,,$$
		and $W$ likewise.
		Similar box bounds can be employed for $c,g$, since $r(x^\star)=0\ll 10^{8}$.
	
	\item can be enforced. Functions that are not Lipschitz continuous, e.g.\ the square-root or Heaviside function, can be made so by replacing them with smoothed functions, e.g.\ via a suitable mollifier. This is a common practice to ensure the derivatives used in a nonlinear optimization algorithm are globally well-defined.
	
	\item can be ensured in practice by using suitable box constraints \eqref{eqn:OCPd}. Note that this assumption also rules out solutions with finite escape times.
	\item is reasonable because without it no numerical method based on piecewise polynomial representation of the solution can converge.
\end{enumerate}
With the above, we give a reduced version of \cite[Thm~2]{PBF}: 
\begin{theorem}
\label{thm:globconv}
	Suppose (A.1)--(A.5) hold.
	Let $\cX_{h,p}$ consist of finite elements of fixed degree $p$ on a quasi-uniform mesh. If $x_{h}$ is a solution to \eqref{eqn:PBF}, then
	\begin{align*}
		g_\text{opt}(x_h),r_\text{feas}(x_h) \xrightarrow{h\rightarrow 0} 0\,.
	\end{align*}
\end{theorem}
\vspace{2mm}

In \cite[Thm~3]{PBF} the result is extended to an order-of-convergence result by extending (A.5) with regard to the order of convergence of the approximation error. The result essentially says: If the finite element approximation converges of order $\ell \in (0.5,p]$ then the optimality gap and feasibility residual converge of the order $\ell/2$.

In addition, \cite{PBF} analyses the convergence when \eqref{eqn:PBF} is solved inexactly, which occurs when integrals are approximated with numerical quadrature.

\section{Implementation
}\label{sec:NLP}
This section outlines how PBF can be implemented.
We first discuss suitable finite element bases. For success of the discretization, a stability condition on the quadrature rule is introduced. Using the basis and quadrature, the NLP is stated and discussed in terms of optimality conditions and sparsity of derivative matrices.

\subsection{Finite Element Basis}
Independent of a collocation method or integral penalty method, a candidate $x_h \in \cX_{h,p}$ is represented using a finite-dimensional solution vector $\bx \in \R^n$, where~$n$ is the total number of finite element functions of all species of $x_h$, which is identified with $\bx$ by means of a basis. Due to considerations on sparsity in the derivative matrices of the NLP, collocation methods use a nodal basis with abscissae that match the collocation points. For an integral penalty method we may opt for a nodal basis with abscissae that are a subset of the quadrature points used for evaluating the integrals; instead, we used equidistant nodes in favour of conditioning in the experiments with PBF in Section~\ref{sec:NumExp}.

\begin{figure}
	\centering
	\includegraphics[width=0.7\columnwidth]{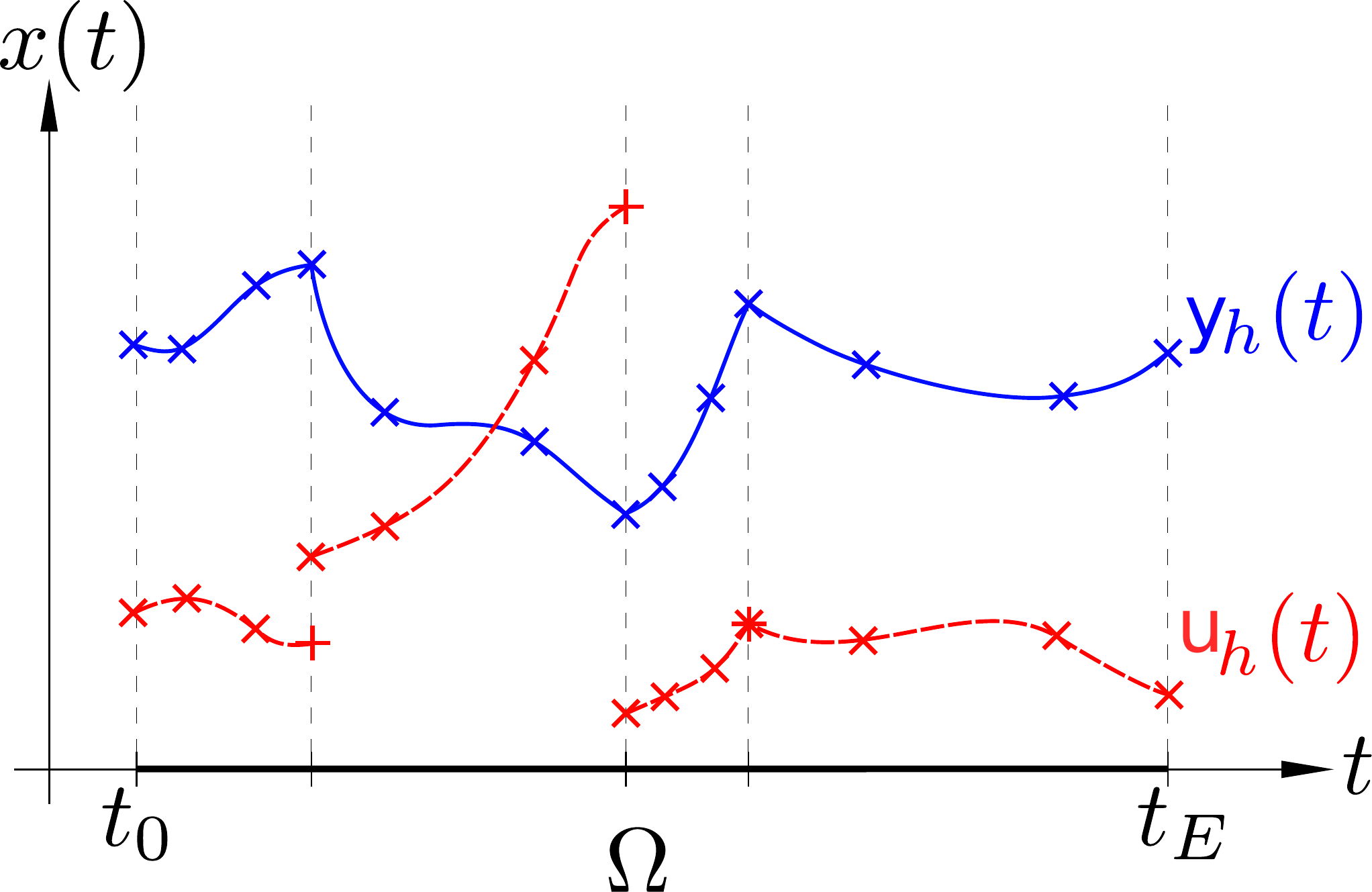}
	\caption{Continuous and discontinuous piecewise polynomial finite element functions $y_h,u_h$ on a mesh of four intervals.}
	\label{fig:finiteelementsyz}
\end{figure}

Species of $y$ must be approximated with continuous finite element functions, whereas for $u$ one can use discontinuous elements. Figure~\ref{fig:finiteelementsyz} shows an example of a nodal basis representation with piecewise polynomial finite element functions of degree $p=3$.

\subsection{Stability Condition and Quadrature Rule}

To minimize \eqref{eqn:PBF}, we approximate the integrals in~$J$ and $r$ with quadrature formulas. As a consequence of this, our numerical solutions $x_h$ will not be exact minimizers of $J_{\omega,\tau}$, but of its numerical quadrature approximation $J_{\omega,\tau,h}$. In order to still be able to prove convergence of~$x_h$, the quadrature rule must be chosen to satisfy the following non-trivial requirement \cite{PBF}:
\begin{align}
	|J_{\omega,\tau}(x_h) - J_{\omega,\tau,h}(x_h)| \leq const \cdot \frac{h^\ell}{\omega}\quad \forall x_h \in \cX_{h,p}\,,\label{eqn:QuadStab}
\end{align}
where $\ell$ is the desired order of convergence of $x_h$. While in general it is not so difficult to devise a convergent quadrature rule, the devil is in the detail that the scheme must converge for \textit{every} element in $\cX_{h,p}$.

\begin{figure}
	\centering
	\includegraphics[width=0.7\columnwidth]{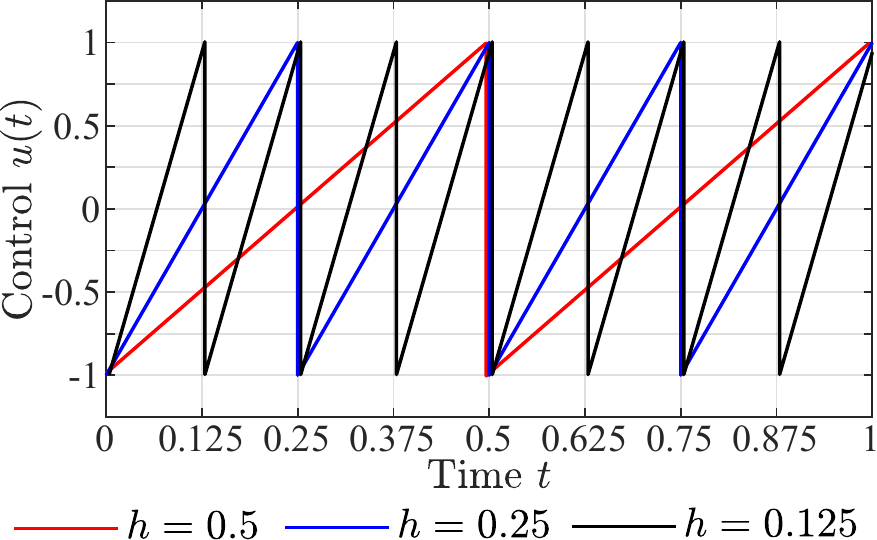}
	\caption{A piecewise affine function, ranging from -1 to 1 on each mesh interval, which serves as a counter-example to satisfaction of \eqref{eqn:QuadStab} when using the midpoint quadrature rule when $p=1$.}
	\label{fig:zigzag}
\end{figure}

We illustrate the difficulty with an example. Consider $\Omega=[0,1]$, and a uniform mesh of size $h \in 1/\N$ of discontinuous piecewise affine finite elements for $u :\Omega\rightarrow \R^1$. Consider the path constraint
$ 	c\big(u(t)\big):=\sin\big(\pi u(t)\big)=0 	$
and approximate $r(x)$ with the midpoint rule quadrature formula as $r_h(x)$. Then \eqref{eqn:QuadStab} requires that the difference
\begin{align*}
	&|r_h(x)-r(x)| \\
	&\ = \Bigg| h \cdot \sum_{j=1}^{1/h} \sin^2\big(\pi u_h(jh-h/2)\big) - \int_{0}^{1}\sin^2\big(\pi u_h(t)\big)\mathrm{d}t \Bigg|
\end{align*}
converges to zero as $h \searrow +0$. However, $\cX_{h,1}$ contains the element $u_h$ plotted in Figure~\ref{fig:zigzag}, which yields $r_h(x)=0$ whereas $r(x)=1/2$. 
Usually, the midpoint quadrature rule converges in $\cO(h^3)$, but since here the integrand $u_h(t)$ can be chosen in a larger space $\cX_{h,1}$, the quadrature rule does not converge as $h$ is decreased. If instead a trapezoidal quadrature rule is used, then \eqref{eqn:QuadStab} is satisfied with order $\ell=3$.

To satisfy \eqref{eqn:QuadStab} in practice, we suggest the use of Gauss-Legendre quadrature with $q=2p$ quadrature points per element. In \cite{PBF} there are further precautions so that, in addition, \eqref{eqn:def:Gamma} can be approximated with Gauss-Legendre quadrature. We write $\alpha_j>0,\tau_j\in\Omega$ for $j=1,\dots,N$ for the quadrature weights and abscissae.

Finally, collocation methods must also satisfy~\eqref{eqn:QuadStab} as well if they want to ensure a bound on the path constraint residual. However, it is not possible for collocation methods to choose $q>p$, due to over-determination of the NLP.
\newcommand{\mE}{{m_{\mathcal{E}}}\xspace}
\newcommand{\mI}{{m_{\mathcal{I}}}\xspace}

\subsection{Nonlinear Program}
We define the following functions and matrix:
\begin{align*}
	F(\bx) {}&{}:= 	W\big(y_h(t_0),y_h(t_F)\big) \\
				&{}\phantom{:=}{}+ \sum_{j=1}^N \alpha_j w\big(y_h(\tau_j),u_h(\tau_j)\big)\,,\tageq\label{eqn:funcF}
				\end{align*}
				\begin{align*}
	C(\bx) {}&{}:= 	\begin{bmatrix}
						g(y_h(t_0),y_h(t_F)\big) \\
						\sqrt{\alpha_1} c\big(\dot{y}_h(\tau_1),y_h(\tau_1),u_h(\tau_1)\big)\\
						\vdots\\
						\sqrt{\alpha_N} c\big(\dot{y}_h(\tau_N),y_h(\tau_N),u_h(\tau_N)\big)
					\end{bmatrix} \in \R^\mE
					\,,\tageq\label{eqn:funcC}
					\end{align*}
				\begin{align*}
	P_s \bx {}&{}:=	\begin{bmatrix}
						s(\tau_1)\\
						\vdots\\
						s(\tau_N)
					\end{bmatrix} \in \R^\mI
					\,,\tageq
\end{align*}
where $\mE:={n_g + q N n_c}$, $\mI:={q N n_s}$.
The quadrature approximation $J_{\omega,\tau,h}$ of the minimization function $J_{\omega,\tau}$ in \eqref{eqn:PBF} can then be expressed as 
\begin{align}
\min_{\bx \in \R^{n}} 	\Phi(\bx):=F(\bx) + \frac{1}{2 \omega} \|C(\bx)\|_2^2 - \tau \bg\t \log(P_s \bx)\,,\tag{NLP}\label{eqn:NLP}
\end{align}
where $\bg \in \R^\mI$ is a vector of quadrature weights.

$\Phi$ should not be minimized with a general-purpose unconstrained minimization solver. Instead, there are suitable tailored solvers for minimization of this exact function with quadratic penalties and log-barriers \cite{PBF,ForsgrenGill,ALM_IPM}. These methods compute vectors \mbox{$\bx \in \R^n,$} \mbox{$\by \in \R^\mE,$} \mbox{$\bz \in \R_+^\mI$,} such that
\begin{align*}
\nabla_\bx \cL(\bx,\by)-P_s\t \bz &= \bO\,,\\
C(\bx) + \omega \by &=\bO\,,\\
\opdiag(P_s \bx) \cdot \bz - \tau \bg &=\bO\,.
\end{align*}
This is equivalent to the first-order necessary optimality conditions for minimizing \eqref{eqn:NLP}, where $\cL(\bx,\by):= F(\bx) - \by\t C(\bx)$.

\subsection{Derivatives and Sparsity}
Second-order NLP solvers need $\nabla_\bx\cL(\bx,\by)$, $\nabla^2_{\bx,\bx}\cL(\bx,\by)$ and $\nabla_\bx C(\bx)$.
We see that \eqref{eqn:funcF}--\eqref{eqn:funcC} essentially use the same formulas as collocation methods, except that they replace collocation points with quadrature points $\tau_j$, and weight the constraints $C$ by $\sqrt{\alpha_j}$. Hence, the Jacobian $\nabla\t_\bx C(\bx) \in \R^{\mE \times n}$ of PBF has the same sparsity structure as for a collocation method of same degree $p$, except the Jacobian in PBF has more rows than in LGR; e.g.\ twice as many rows as a collocation method of the same degree if we choose to use Gauss-Legendre quadrature with $q=2p$ abscissae per element, as proposed.
The observations for $\nabla\t_\bx C(\bx)$ can be extended to $\nabla^2_{\bx,\bx}\cL(\bx,\by)\in \R^{n \times n}$: PBF($p$) and LGR($p$) of the same degree $p$ on the same mesh have the same Hessian dimension and sparsity structure. Section~\ref{sec:NumExp} presents sparsity plots for PBF(5) and LGR(5).

\section{Numerical Experiments}\label{sec:NumExp}
\begin{figure}
	\centering
	\includegraphics[width=0.8\columnwidth]{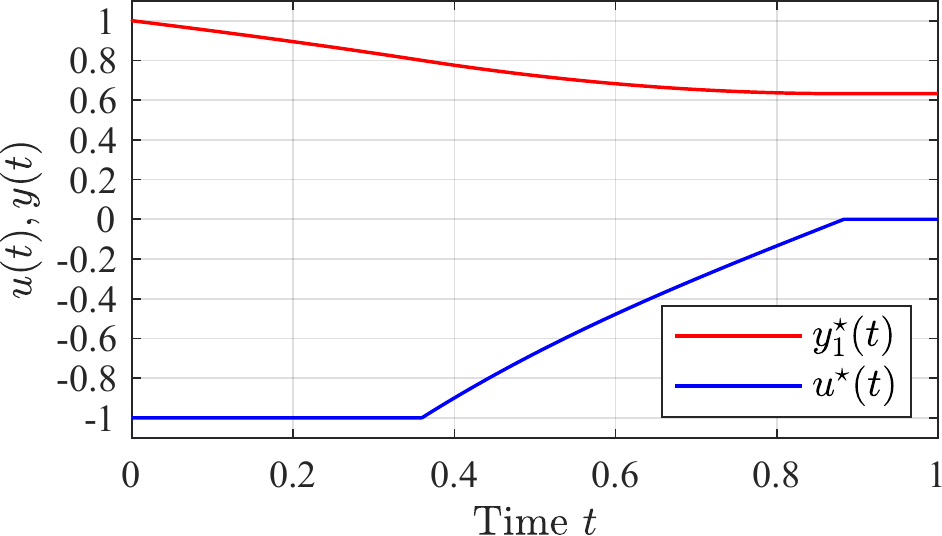}
	\caption{Analytical solution to \eqref{eqn:ExampleOCP2}.}
	\label{fig:exp2}
\end{figure}

We now consider a test problem for which both types of methods converge with success, so that we can compare conditioning, convergence, and rate of convergence to a known analytical solution:
\begin{align*}
	\min_{y,u}& & J&= y_2(1),\tag{X2}\label{eqn:ExampleOCP2}\\
	\text{s.t.}& & y_1(0)&=1,\ \, \dot{y}_1(t)=\frac{u(t)}{2y_1(t)},\ \, \sqrt{0.4}\leq y_1(t),\\
	& & y_2(0)&=0,\ \, \dot{y}_2(t)=4 y_1(t)^4 + u(t)^2,\ \,  -1\leq u(t)\,.
\end{align*}
The solution is shown in \abbfig~\ref{fig:exp2}.
$u^\star$ is constant outside $t_0=1-\frac{\sqrt{41}}{10}\approx 0.35$ and $t_1=t_0+\log 2-\frac{\log(\sqrt{41}-5)}{2}\approx 0.88$, between which $u^\star(t)= 0.8 \sinh\left(2 (t - t_1)\right)$, yielding $J\approx 2.057866062168$. 

\begin{figure}
	\centering
	\includegraphics[width=0.85\columnwidth]{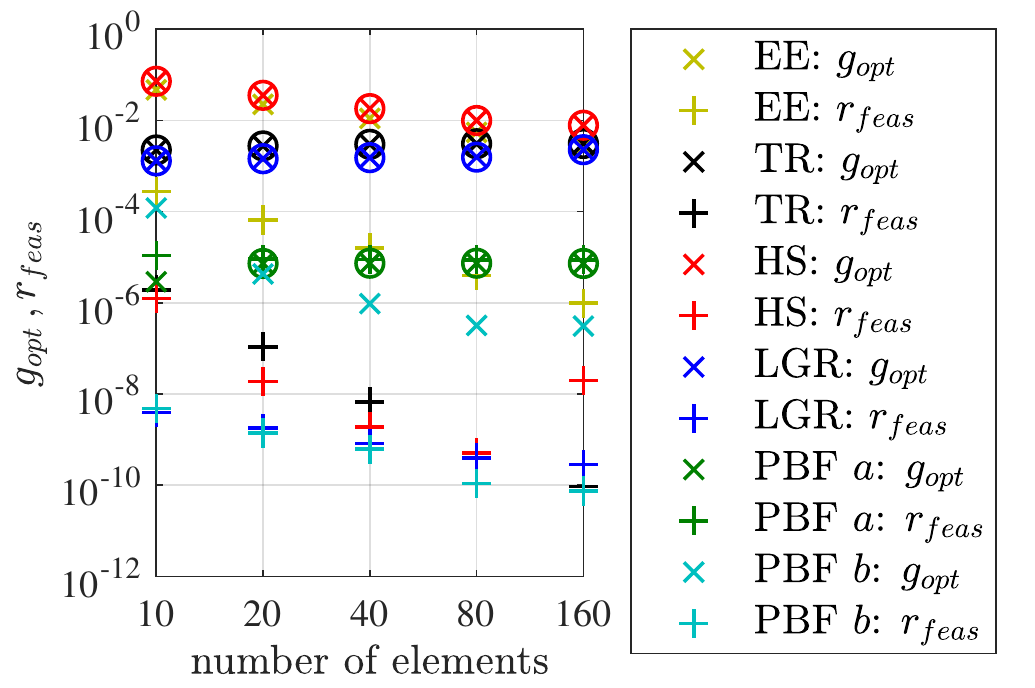}
	\caption{Convergence of optimality gap and feasibility residual
	LGR and PBF use polynomial degree $p=5$. PBF uses two different values: a) $\omega=10^{-5}$, b) $\omega=10^{-10}$.}
	\label{fig:loglogxpl2}
\end{figure}

All methods yield accurate solutions. \abbFig~\ref{fig:loglogxpl2} shows the convergence of the optimality gap and feasibility residual of a respective method. Remarking on the former, we computed $J^\star-J(x_h)$ and encircled the cross when $J(x_h)<J^\star$. Note in the figure that for $\geq 40$ elements the most accurate solutions in terms of feasibility are found by PBF with $\omega=10^{-10}$. Further, we find that the collocation methods significantly underestimate the optimality value for this experiment.

Now we discuss rates of convergence. Convergence of only first order is expected because higher derivatives of $y^\star$ are non-smooth and $u^\star$ has edges. Indeed, $r_\text{feas}$ converges linearly for all methods. PBF5 with $\omega=10^{-5}$ stagnates early because it converges to the optimal penalty solution, which for this instance is converged from $20$ elements onwards. $g_\text{opt},r_\text{feas}$ are then fully determined by $\omega$. The issue is resolved by choosing $\omega$ smaller. LGR5 and PBF5 with $\omega=10^{-10}$ converge similarly, and stagnate at $r_\text{feas}\approx 10^{-10}$. Due to the high exponent in the objective, small feasibility errors in the collocation methods amount to significant underestimation of the objective.

\begin{figure}
	\centering
	\includegraphics[width=0.7\columnwidth]{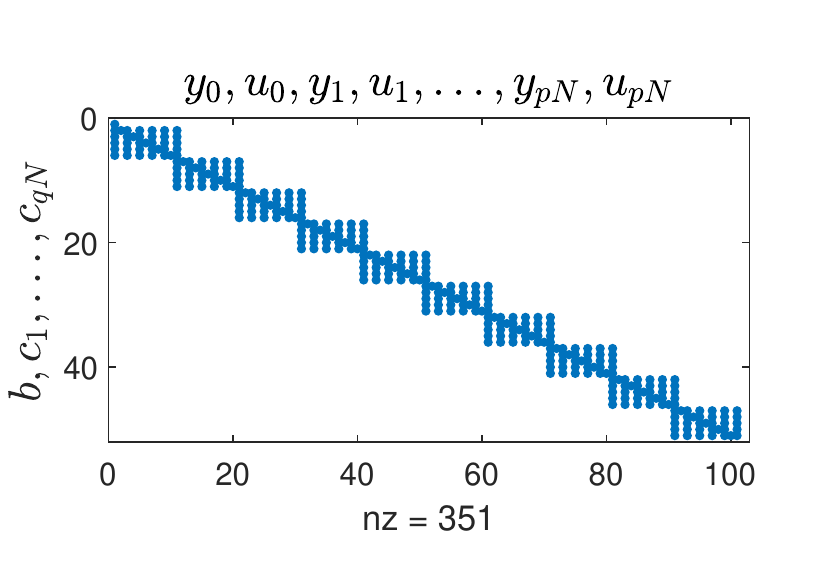}
	\caption{Sparsity of $\nabla_\bx C(\bx)\t$ for LGR5 when $h=\frac{1}{10}$, i.e.~$N=10$. For LGR, notice $q=p-1$. The discretization does not depend on $u_{pN}$.}
	\label{fig:sparsitylgr}
\end{figure}

\begin{figure}
	\centering
	\includegraphics[width=0.7\columnwidth]{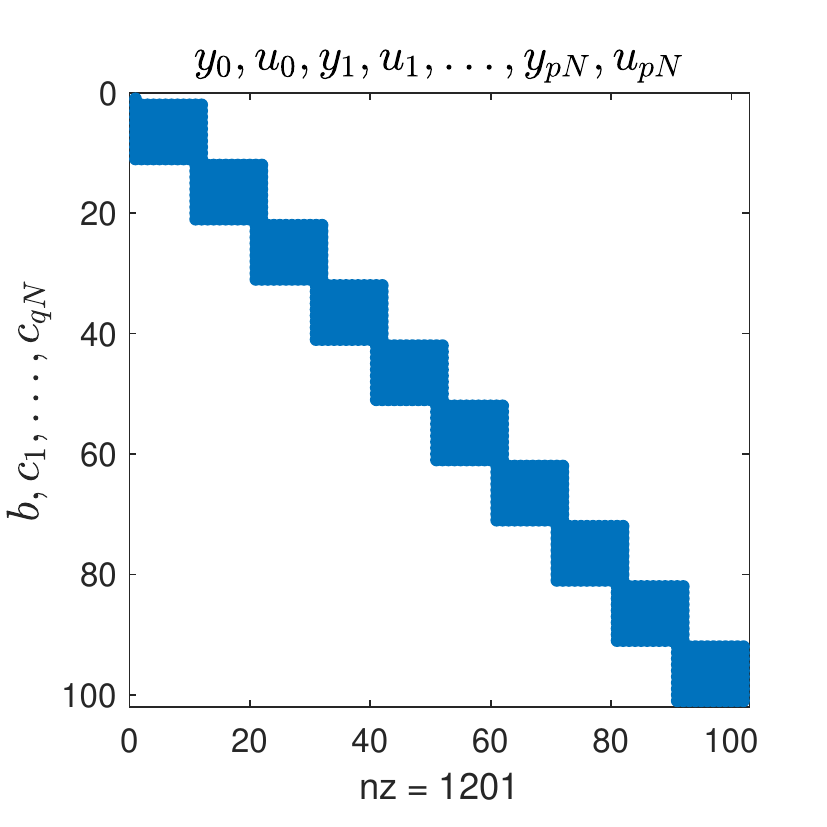}
	\caption{Sparsity of $\nabla_\bx C(\bx)\t$ for PBF5 when $h=\frac{1}{10}$, i.e.~$N=10$, with $q=2p$.}
	\label{fig:sparsitypbf}
\end{figure}

Finally, we look into computational cost. Solving the collocation methods with IPOPT and the PBF5 discretization with the interior-point method in \cite{ForsgrenGill}, the optimization converges in $\approx 20$ iterations for any discretization. Differences in computational cost can arise when one discretization results in much denser or larger problems than others. Here, we compare the sparsity structure of the Jacobian $\nabla_\bx C(\bx)\t$ for LGR5 in \abbfig~\ref{fig:sparsitylgr} and PBF5 in \abbfig~\ref{fig:sparsitypbf}, each using a mesh size of $h=\frac{1}{10}$.
Note that for PBF5, $C(\bx)$ has a larger dimension $\mE$ than LGR5, thus the Jacobian has more rows and hence more non-zeros. However, critical for computations is the primal Schur complement $\mathbf{\Sigma}=\nabla_{\bx\bx}^2 \cL(\bx,\blambda)+\nabla_\bx C(\bx)\t \bD \nabla_\bx C(\bx)$, which is used when solving the KKT system via the reduced form, where $\bD$ is a diagonal matrix. $\mathbf{\Sigma}$ is a narrow-banded matrix with dense band of the same bandwidth for LGR5 and PBF5.

With regard to computational cost, it follows from \abbfig~\ref{fig:loglogxpl2} that the ability to choose $\omega$ in PBF can be advantageous. In particular, on coarse meshes, one may opt for small feasibility residual by manually decreasing $\omega$, whereas with a collocation method one is stuck with the feasibility residual that one obtains for that particular mesh. The figure shows this: For $\omega=10^{-10}$, even on the coarsest mesh the PBF method achieves a solution that has a smaller feasibility residual than other methods on the same mesh. For this problem this becomes possible because the path constraint could be satisfied with zero error by choosing $y$ a polynomial of degree 3 (because here PBF uses $p=5$).

\section{Conclusions}
\balance	
We compared classical direct discretization methods, i.e.\ collocation-type methods, with integral penalty-type methods, and in particular with PBF. We highlighted the critical difference, that is, forcing exact constraint satisfaction at a finite number of points for the former, versus forcing approximate constraint satisfaction everywhere for the latter. We also described the consequences of these approaches on the convergence: Collocation-type methods achieve a fixed bias between optimality gap and feasibility residual at a given mesh, whereas with PBF one can trade between the two. Further, collocation methods may fail to converge in several cases where PBF converges.

The convergence has been investigated by reviewing theorems and performing numerical experiments. Both indicate that direct transcription based on integral penalties is competitive in computational cost and accuracy compared to state-of-the-art collocation methods in the general case, but in addition also works reliably in applications where the latter struggle to converge.

The robustness of an optimal control solver is determined by the weakest link in the chain. The discretization is only one half of the numerical method. Future work could therefore be dedicated to the development of tailored NLP solvers to improve reliability and computational efficiency.


\bibliographystyle{IEEEtran} 	
\bibliography{PBFconf_refs_Martin_Eric}

\begin{thebibliography}{10}
\providecommand{\url}[1]{#1}
\csname url@rmstyle\endcsname
\providecommand{\newblock}{\relax}
\providecommand{\bibinfo}[2]{#2}
\providecommand\BIBentrySTDinterwordspacing{\spaceskip=0pt\relax}
\providecommand\BIBentryALTinterwordstretchfactor{4}
\providecommand\BIBentryALTinterwordspacing{\spaceskip=\fontdimen2\font plus
\BIBentryALTinterwordstretchfactor\fontdimen3\font minus
  \fontdimen4\font\relax}
\providecommand\BIBforeignlanguage[2]{{%
\expandafter\ifx\csname l@#1\endcsname\relax
\typeout{** WARNING: IEEEtran.bst: No hyphenation pattern has been}%
\typeout{** loaded for the language `#1'. Using the pattern for}%
\typeout{** the default language instead.}%
\else
\language=\csname l@#1\endcsname
\fi
#2}}

\bibitem{Betts2nd}
J.~T. Betts, \emph{Practical Methods for Optimal Control and Estimation Using
  Nonlinear Programming}, 2nd~ed.\hskip 1em plus 0.5em minus 0.4em\relax New
  York, NY, USA: Cambridge University Press, 2010.

\bibitem{Balakrishnan68}
A.~V. Balakrishnan, ``On a new computing technique in optimal control,''
  \emph{SIAM J. Control}, vol.~6, pp. 149--173, 1968.

\bibitem{Hager90}
\BIBentryALTinterwordspacing
W.~W. Hager, ``Multiplier methods for nonlinear optimal control,'' \emph{SIAM
  J. Numer. Anal.}, vol.~27, no.~4, pp. 1061--1080, 1990. [Online]. Available:
  \url{https://doi.org/10.1137/0727063}
\BIBentrySTDinterwordspacing

\bibitem{PBF}
M.~P. Neuenhofen and E.~C. Kerrigan, ``Dynamic optimization with convergence
  guarantees,'' \emph{arXiv:1810.04059}, 2018.

\bibitem{KellyMatthew}
\BIBentryALTinterwordspacing
M.~Kelly, ``An introduction to trajectory optimization: How to do your own
  direct collocation,'' \emph{SIAM Rev.}, vol.~59, no.~4, pp. 849--904, 2017.
  [Online]. Available: \url{https://doi.org/10.1137/16M1062569}
\BIBentrySTDinterwordspacing

\bibitem{GPOPS2}
\BIBentryALTinterwordspacing
M.~A. Patterson and A.~V. Rao, ``{{GPOPS}-{II}}: a {M}atlab software for
  solving multiple-phase optimal control problems using {$hp$}-adaptive
  {G}aussian quadrature collocation methods and sparse nonlinear programming,''
  \emph{ACM Trans. Math. Software}, vol.~41, no.~1, pp. Art. 1, 37, 2014.
  [Online]. Available: \url{https://doi.org/10.1145/2558904}
\BIBentrySTDinterwordspacing

\bibitem{RaoHager}
\BIBentryALTinterwordspacing
W.~W. Hager, J.~Liu, S.~Mohapatra, A.~V. Rao, and X.-S. Wang, ``Convergence
  rate for a {G}auss collocation method applied to constrained optimal
  control,'' \emph{SIAM J. Control Optim.}, vol.~56, no.~2, pp. 1386--1411,
  2018. [Online]. Available: \url{https://doi.org/10.1137/16M1096761}
\BIBentrySTDinterwordspacing

\bibitem{Conway2012}
\BIBentryALTinterwordspacing
B.~A. Conway, ``{A Survey of Methods Available for the Numerical Optimization
  of Continuous Dynamic Systems},'' \emph{Journal of Optimization Theory and
  Applications}, vol. 152, no.~2, pp. 271--306, Feb 2012. [Online]. Available:
  \url{https://doi.org/10.1007/s10957-011-9918-z}
\BIBentrySTDinterwordspacing

\bibitem{rao_asurvey}
A.~V. Rao, ``{Survey of Numerical Methods for Optimal Control},'' in
  \emph{{Advances in the Astronautical Sciences}}, vol. 135, 2010.

\bibitem{Quirynen2017}
R.~Quirynen, S.~Gros, B.~Houska, and M.~Diehl, ``Lifted collocation integrators
  for direct optimal control in acado toolkit,'' \emph{Mathematical Programming
  Computation}, vol.~9, no.~4, p. 527–571, 2017.

\bibitem{Maurer}
K.~Malanowski, C.~B\"{u}skens, and H.~Maurer, ``Convergence of approximations
  to nonlinear optimal control problems,'' in \emph{Mathematical programming
  with data perturbations}, ser. Lecture Notes in Pure and Appl. Math.\hskip
  1em plus 0.5em minus 0.4em\relax Dekker, New York, 1998, vol. 195, pp.
  253--284.

\bibitem{kameswaran_biegler_2008}
S.~Kameswaran and L.~T. Biegler, ``Convergence rates for direct transcription
  of optimal control problems using collocation at {R}adau points,''
  \emph{Comput. Optim. Appl.}, vol.~41, pp. 81--126, 2008.

\bibitem{ForsgrenGill}
\BIBentryALTinterwordspacing
A.~Forsgren and P.~E. Gill, ``Primal-dual interior methods for nonconvex
  nonlinear programming,'' \emph{SIOPT}, vol.~8, no.~4, pp. 1132--1152, 1998.
  [Online]. Available: \url{https://doi.org/10.1137/S1052623496305560}
\BIBentrySTDinterwordspacing

\bibitem{ALM_IPM}
\BIBentryALTinterwordspacing
Y.~Cao, A.~Seth, and C.~D. Laird, ``An augmented {L}agrangian interior-point
  approach for large-scale {NLP} problems on graphics processing units,''
  \emph{Computers and Chemical Engineering}, vol.~85, pp. 76 -- 83, 2016.
  [Online]. Available:
  \url{http://www.sciencedirect.com/science/article/pii/S0098135415003257}
\BIBentrySTDinterwordspacing

\end{thebibliography}

\end{document}